\documentclass[11pt]{article}
\usepackage{amsfonts}
\usepackage{latexsym}
\usepackage{cite}
\usepackage{amsmath,amsfonts,latexsym,amssymb}
\usepackage[mathscr]{eucal}
\usepackage{cases,color}
\usepackage{amsthm}

\usepackage[bf,small]{caption2}
\usepackage{float}
\usepackage{graphicx}
\usepackage{amsmath}
\usepackage{amssymb}
\usepackage[all]{xy}

\newtheorem{theorem}{theorem}[section]

\newtheorem{thm}[theorem]{Theorem}

\begin{document}

\title{\textbf{The symmetric slice of ${\rm SL}(3,\mathbb{C})$-character variety of the Whitehead link}}
\author{\Large Haimiao Chen}
\date{}
\maketitle

\begin{abstract}
  We give a nice description for a Zariski open subset of the ${\rm SL}(3,\mathbb{C})$-character variety of the Whitehead link.

\end{abstract}

\section{Introduction}

For a finitely presented group $\Gamma$, let $\mathcal{X}_3(\Gamma)$ denote the irreducible ${\rm SL}(3,\mathbb{C})$-character variety of $\Gamma$, which parameterizes conjugacy classes of irreducible representations $\Gamma\to{\rm SL}(3,\mathbb{C})$.

While various results on ${\rm SL}(2,\mathbb{C})$-character varieties can be seen in the literature, the results for
${\rm SL}(3,\mathbb{C})$-character varieties are quite rare (see \cite{La07,HMP16,MP16}, etc). It is now still considered as a difficult task to determine $\mathcal{X}_3(\Gamma)$, even when $\Gamma$ is generated by two elements.

From now on, let $\Gamma$ be the fundamental group of the exterior of the Whitehead link.

\begin{figure}[H]
  \centering
  \includegraphics[width=9cm]{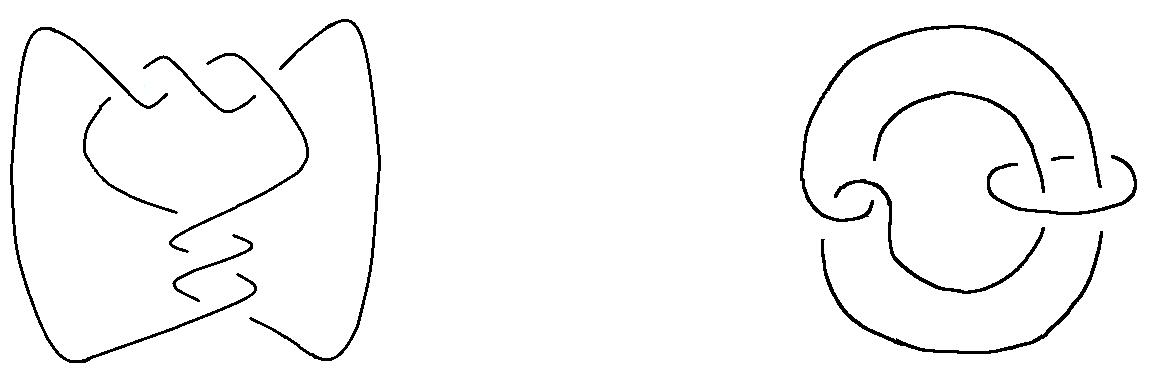}\\
  \caption{Two projection diagrams of the Whitehead link.}\label{fig:Whitehead}
\end{figure}

Following \cite{Ch22-1} Section 4, the left projection diagram in Figure \ref{fig:Whitehead} gives rise to a presentation
$\Gamma\cong\langle y,z\mid yzy\overline{z}^2yz=zy\overline{z}^2yzy\rangle.$

We use boldface letters to denote elements of ${\rm SL}(3,\mathbb{C})$. Let $\mathbf{e}$ denote the identity matrix. Given $\mathbf{x}\in{\rm SL}(3,\mathbb{C})$, let $\mathbf{x}^{\rm tr}$ denote its transpose.

Given $\mathbf{y},\mathbf{z}\in{\rm SL}(3,\mathbb{C})$, a sufficient and necessary condition for there to exist a representation $\rho:\pi\to{\rm SL}(3,\mathbb{C})$ with $\rho(y)=\mathbf{y}$, $\rho(z)=\mathbf{z}$ is
\begin{align}
\mathbf{y}\mathbf{z}\mathbf{y}\overline{\mathbf{z}}^2\mathbf{y}\mathbf{z}=\mathbf{z}\mathbf{y}\overline{\mathbf{z}}^2\mathbf{y}\mathbf{z}\mathbf{y}.
\label{eq:relation-0}
\end{align}
When this holds, $\rho$ is determined by $\mathbf{y},\mathbf{z}$. The irreducibility of $\rho$ is equivalent to that $(\mathbf{y},\mathbf{z})$ is irreducible, i.e. neither $\mathbf{y},\mathbf{z}$ nor $\mathbf{y}^{\rm tr},\mathbf{z}^{\rm tr}$ share an eigenvector.

Call the subset $\mathcal{X}_3^{\rm asym}(\Gamma)$ (resp. $\mathcal{X}_3^{\rm sym}(\Gamma)$) contributed by $\mathbf{y},\mathbf{z}$ with
${\rm tr}[\mathbf{y},\mathbf{z}]\ne{\rm tr}([\mathbf{z},\mathbf{y}])$ (resp. ${\rm tr}[\mathbf{y},\mathbf{z}]={\rm tr}([\mathbf{z},\mathbf{y}])$) the {\it asymmetric slice} (resp. {\it symmetric slice}) of $\mathcal{X}_3(\Gamma)$. A nice description for $\mathcal{X}_3^{\rm asym}(\Gamma)$ had been given in \cite{GW19}. In \cite{Ch22-1} we recovered $\mathcal{X}_3^{\rm asym}(\Gamma)$, and derived equations defining $\mathcal{X}_3^{\rm sym}(\Gamma)$ (see Theorem 4.3 there). However, those equations are too complicated to yield a good understanding of
$\mathcal{X}_3^{\rm sym}(\Gamma)$.

As the main result of this paper, we show
\begin{thm}\label{thm:main}
The symmetric slice of the ${\rm SL}(3,\mathbb{C})$-character variety of the Whitehead link contains a Zariski open subset which is a regular $\mathbb{Z}/6\mathbb{Z}$-cover of a Zariski open subset of some hypersurface in $\mathbb{C}^5$.
\end{thm}




\section{A geometric insight towards the symmetric slice}

\subsection{Preparation}

Let $\mathcal{M}^0$, $\mathcal{M}^1$ respectively denote the vector space of symmetric and skew-symmetric $3\times 3$ matrices.

Identify $\mathbf{u}=(u_{ij})_{3\times 3}\in\mathcal{M}^1$ with $\mathbf{u}^\neg:=(u_{12},u_{13},u_{23})^{\rm tr}\in\mathbb{C}^3$.

It is easy to see that for any $\mathbf{u},\mathbf{v},\mathbf{w}\in\mathcal{M}^1$,
\begin{align*}
(\mathbf{u}\mathbf{v}-\mathbf{v}\mathbf{u})^\neg&=\mathbf{u}^\neg\times\mathbf{v}^\neg,  \\
{\rm tr}(\mathbf{u}\mathbf{v}\mathbf{w})&=\det[\mathbf{u}^\neg,\mathbf{v}^\neg,\mathbf{w}^\neg],
\end{align*}
Consequently, $\mathbf{u},\mathbf{v}$ are linearly dependent if and only if $\mathbf{u}\mathbf{v}=\mathbf{v}\mathbf{u}$, and $\mathbf{u},\mathbf{v},\mathbf{w}$ are linearly dependent if and only if ${\rm tr}(\mathbf{u}\mathbf{v}\mathbf{w})=0$.

Call $\mathbf{x}$ {\it ordinary} if its minimal polynomial is $\det(\lambda\mathbf{e}-\mathbf{x})$, or equivalently, each eigen-subspace is $1$-dimensional; otherwise, call $\mathbf{x}$ {\it special}.

For any $\mathbf{a}\in{\rm SL}(3,\mathbb{C})$, by Hamilton-Cayley Theorem,
\begin{align}
\mathbf{a}^2&={\rm tr}(\mathbf{a})\cdot\mathbf{a}-{\rm tr}(\overline{\mathbf{a}})\cdot\mathbf{e}+\overline{\mathbf{a}}, \label{eq:square}  \\
\mathbf{a}^3&=({\rm tr}(\mathbf{a})^2-{\rm tr}(\overline{\mathbf{a}}))\cdot\mathbf{a}
+(1-{\rm tr}(\mathbf{a}){\rm tr}(\overline{\mathbf{a}}))\cdot\mathbf{e}+{\rm tr}(\mathbf{a})\cdot\overline{\mathbf{a}}.  \label{eq:cube}
\end{align}
For any $\mathbf{a},\mathbf{b}\in{\rm SL}(3,\mathbb{C})$, we have
\begin{align}
\mathbf{a}\mathbf{b}\mathbf{a}&=-\overline{\mathbf{a}}\mathbf{b}-\mathbf{b}\overline{\mathbf{a}}+{\rm tr}(\mathbf{a}\mathbf{b})\cdot\mathbf{a}
+{\rm tr}(\mathbf{b})\cdot\overline{\mathbf{a}}+{\rm tr}(\overline{\mathbf{a}})\cdot\mathbf{b}  \nonumber  \\
&\ \ \ \ +({\rm tr}(\overline{\mathbf{a}}\mathbf{b})-{\rm tr}(\overline{\mathbf{a}}){\rm tr}(\mathbf{b}))\cdot\mathbf{e}, \label{eq:identity-1}  \\
\mathbf{a}\mathbf{b}\mathbf{a}&={\rm tr}(\mathbf{a}\mathbf{b})\cdot\mathbf{a}
-{\rm tr}(\overline{\mathbf{a}}\overline{\mathbf{b}})\cdot\overline{\mathbf{b}}+\overline{\mathbf{b}}\overline{\mathbf{a}}\overline{\mathbf{b}}.
\label{eq:identity-2}
\end{align}
The first identity is simplified from (6) in \cite{Ch22-1}, by applying (\ref{eq:square}).

\subsection{A Zariski open subset of $\mathcal{X}_3^{\rm sym}(\Gamma)$}

By \cite{Ch22-1} Theorem 3.8 and Theorem 3.9, we may just assume $\mathbf{y},\mathbf{z}\in\mathcal{M}^0$.

Let $\mathbf{a}=\overline{\mathbf{y}}\mathbf{z}$, $\mathbf{b}=\mathbf{z}\overline{\mathbf{y}}$.
Then $\mathbf{z}=\mathbf{y}\mathbf{a}=\mathbf{b}\mathbf{y}$ and (\ref{eq:relation-0}) becomes
$$\mathbf{y}[\mathbf{a},\overline{\mathbf{b}}]=[\mathbf{a},\overline{\mathbf{b}}]\mathbf{y}.$$

Note that $\mathbf{b}=\mathbf{a}^{\rm tr}$. Let
\begin{align*}
t={\rm tr}(\mathbf{a})={\rm tr}(\mathbf{b}), \qquad \overline{t}={\rm tr}(\overline{\mathbf{a}})={\rm tr}(\overline{\mathbf{b}}), \\
r={\rm tr}(\overline{\mathbf{a}}\mathbf{b})={\rm tr}(\mathbf{a}\overline{\mathbf{b}}),  \qquad
s={\rm tr}(\mathbf{a}\mathbf{b}), \qquad  \overline{s}={\rm tr}(\overline{\mathbf{a}}\overline{\mathbf{b}}).
\end{align*}
Denote $t_{1^2\overline{2}}={\rm tr}(\mathbf{a}^2\overline{\mathbf{b}})$,
$t_{1\overline{2}\overline{1}^22}={\rm tr}(\mathbf{a}\overline{\mathbf{b}}\overline{\mathbf{a}}^2\mathbf{b})$, and so forth.

Assume $[\mathbf{a},\overline{\mathbf{b}}]$ is ordinary. Since $\mathbf{y}$ commutes with $[\mathbf{a},\overline{\mathbf{b}}]$, by \cite{Ch22-1} Lemma 2.2 we may write $$\mathbf{y}=\lambda[\mathbf{a},\overline{\mathbf{b}}]+\mu\mathbf{e}+\nu[\overline{\mathbf{b}},\mathbf{a}]$$
for some $\lambda,\mu,\nu\in\mathbb{C}$.
Then $\mathbf{y}\mathbf{a}=\mathbf{b}\mathbf{y}$ becomes
\begin{align}
\lambda(\mathbf{a}\overline{\mathbf{b}}\overline{\mathbf{a}}\mathbf{b}\mathbf{a}
-\mathbf{b}\mathbf{a}\overline{\mathbf{b}}\overline{\mathbf{a}}\mathbf{b})
+\mu(\mathbf{a}-\mathbf{b})+\nu(\overline{\mathbf{b}}\mathbf{a}\mathbf{b}-\mathbf{a}\mathbf{b}\overline{\mathbf{a}})=0.   \label{eq:equaltiy-1}
\end{align}

Observe that $\mathbf{a}\overline{\mathbf{b}}\overline{\mathbf{a}}\mathbf{b}\mathbf{a}
-\mathbf{b}\mathbf{a}\overline{\mathbf{b}}\overline{\mathbf{a}}\mathbf{b}, \mathbf{a}-\mathbf{b},
\overline{\mathbf{b}}\mathbf{a}\mathbf{b}-\mathbf{a}\mathbf{b}\overline{\mathbf{a}}\in\mathcal{M}^1$,
so the equation (\ref{eq:equaltiy-1}) has a nontrivial solution if and only if
\begin{align*}
0&={\rm tr}\big((\mathbf{a}\overline{\mathbf{b}}\overline{\mathbf{a}}\mathbf{b}\mathbf{a}
-\mathbf{b}\mathbf{a}\overline{\mathbf{b}}\overline{\mathbf{a}}\mathbf{b})(\mathbf{a}-\mathbf{b})
(\overline{\mathbf{b}}\mathbf{a}\mathbf{b}-\mathbf{a}\mathbf{b}\overline{\mathbf{a}})\big)   \\
&={\rm tr}\big(\mathbf{a}\overline{\mathbf{b}}\overline{\mathbf{a}}\mathbf{b}\mathbf{a}
\big((\mathbf{a}-\mathbf{b})(\overline{\mathbf{b}}\mathbf{a}\mathbf{b}-\mathbf{a}\mathbf{b}\overline{\mathbf{a}})
-(\overline{\mathbf{b}}\mathbf{a}\mathbf{b}-\mathbf{a}\mathbf{b}\overline{\mathbf{a}})(\mathbf{a}-\mathbf{b})\big)\big)  \\
&={\rm tr}\big(\mathbf{a}\overline{\mathbf{b}}\overline{\mathbf{a}}\mathbf{b}\mathbf{a}(\mathbf{a}\overline{\mathbf{b}}\mathbf{a}\mathbf{b}
-\mathbf{a}\mathbf{b}\overline{\mathbf{a}}\mathbf{b}+\overline{\mathbf{b}}\mathbf{a}\mathbf{b}^2-\mathbf{a}^2\mathbf{b}\overline{\mathbf{a}}
+\mathbf{b}\mathbf{a}\mathbf{b}\overline{\mathbf{a}}-\overline{\mathbf{b}}\mathbf{a}\mathbf{b}\mathbf{a})\big)=:K.
\end{align*}

Luckily,
\begin{align*}
{\rm tr}(\mathbf{a}\overline{\mathbf{b}}\overline{\mathbf{a}}\mathbf{b}\mathbf{a}\cdot\mathbf{a}\mathbf{b}\overline{\mathbf{a}}\mathbf{b})
&={\rm tr}((\mathbf{a}\overline{\mathbf{b}}\overline{\mathbf{a}}\mathbf{b}\mathbf{a}^2\mathbf{b}\overline{\mathbf{a}}\mathbf{b})^{\rm tr})
={\rm tr}(\mathbf{a}\overline{\mathbf{b}}\mathbf{a}\mathbf{b}^2\mathbf{a}\overline{\mathbf{b}}\overline{\mathbf{a}}\mathbf{b})
={\rm tr}(\mathbf{a}\overline{\mathbf{b}}\overline{\mathbf{a}}\mathbf{b}\mathbf{a}\cdot\overline{\mathbf{b}}\mathbf{a}\mathbf{b}^2),  \\
{\rm tr}(\mathbf{a}\overline{\mathbf{b}}\overline{\mathbf{a}}\mathbf{b}\mathbf{a}\cdot\mathbf{a}^2\mathbf{b}\overline{\mathbf{a}})
&={\rm tr}(\mathbf{a}^2\mathbf{b})={\rm tr}((\mathbf{a}^2\mathbf{b})^{\rm tr})={\rm tr}(\mathbf{a}\mathbf{b}^2)
={\rm tr}(\mathbf{a}\overline{\mathbf{b}}\overline{\mathbf{a}}\mathbf{b}\mathbf{a}\cdot\mathbf{b}\mathbf{a}\mathbf{b}\overline{\mathbf{a}}).
\end{align*}
Hence $K$ equals
\begin{align*}
&{\rm tr}(\mathbf{a}\overline{\mathbf{b}}\overline{\mathbf{a}}\mathbf{b}\mathbf{a}\cdot\mathbf{a}\overline{\mathbf{b}}\mathbf{a}\mathbf{b})
-{\rm tr}(\mathbf{a}\overline{\mathbf{b}}\overline{\mathbf{a}}\mathbf{b}\mathbf{a}\cdot\overline{\mathbf{b}}\mathbf{a}\mathbf{b}\mathbf{a})
={\rm tr}(\mathbf{a}\overline{\mathbf{b}}\overline{\mathbf{a}}\mathbf{b}\mathbf{a}^2\overline{\mathbf{b}}\mathbf{a}\mathbf{b})
-{\rm tr}(\mathbf{a}^2\overline{\mathbf{b}}\overline{\mathbf{a}}\mathbf{b}\mathbf{a}\overline{\mathbf{b}}\mathbf{a}\mathbf{b})   \\
=\ &{\rm tr}\big(\mathbf{a}\overline{\mathbf{b}}\overline{\mathbf{a}}\mathbf{b}(t\mathbf{a}-\overline{t}\mathbf{e}+\overline{\mathbf{a}})
\overline{\mathbf{b}}\mathbf{a}\mathbf{b}\big)
-{\rm tr}\big((t\mathbf{a}-\overline{t}\mathbf{e}+\overline{\mathbf{a}})
\overline{\mathbf{b}}\overline{\mathbf{a}}\mathbf{b}\mathbf{a}\overline{\mathbf{b}}\mathbf{a}\mathbf{b}\big)  \\
=\ &{\rm tr}(\mathbf{a}\overline{\mathbf{b}}\overline{\mathbf{a}}\mathbf{b}\overline{\mathbf{a}}\overline{\mathbf{b}}\mathbf{a}\mathbf{b})
-{\rm tr}(\overline{\mathbf{a}}\overline{\mathbf{b}}\overline{\mathbf{a}}\mathbf{b}\mathbf{a}\overline{\mathbf{b}}\mathbf{a}\mathbf{b})  \\
=\ &{\rm tr}(\mathbf{a}\overline{\mathbf{b}}(r\overline{\mathbf{a}}-r\overline{\mathbf{b}}
+\overline{\mathbf{b}}\mathbf{a}\overline{\mathbf{b}})\overline{\mathbf{b}}\mathbf{a}\mathbf{b})
-{\rm tr}(\overline{\mathbf{a}}\overline{\mathbf{b}}\overline{\mathbf{a}}\mathbf{b}(r\mathbf{a}-r\mathbf{b}
+\mathbf{b}\overline{\mathbf{a}}\mathbf{b})\mathbf{b})   \\
=\ &-r\cdot{\rm tr}(\mathbf{a}\overline{\mathbf{b}}^3\mathbf{a}\mathbf{b})+{\rm tr}((\mathbf{a}\overline{\mathbf{b}}^2)^2\mathbf{a}\mathbf{b})
+r\cdot{\rm tr}(\overline{\mathbf{a}}\overline{\mathbf{b}}\overline{\mathbf{a}}\mathbf{b}^3)
-{\rm tr}(\overline{\mathbf{a}}\overline{\mathbf{b}}(\overline{\mathbf{a}}\mathbf{b}^2)^2)  \\
=\ &-r\big((\overline{t}^2-t){\rm tr}(\mathbf{a}\mathbf{b}\mathbf{a}\overline{\mathbf{b}})+(1-t\overline{t}){\rm tr}(\mathbf{a}^2\mathbf{b})
+\overline{t}\cdot{\rm tr}((\mathbf{a}\mathbf{b})^2)\big)   \\
&+r\big((t^2-\overline{t}){\rm tr}(\overline{\mathbf{a}}\overline{\mathbf{b}}\overline{\mathbf{a}}\mathbf{b})
+(1-t\overline{t}){\rm tr}(\overline{\mathbf{a}}^2\overline{\mathbf{b}})+t\cdot{\rm tr}((\overline{\mathbf{a}}\overline{\mathbf{b}})^2)\big)  \\
&+{\rm tr}((t_{1\overline{2}^2}\mathbf{a}\overline{\mathbf{b}}^2-t_{\overline{1}2^2}\mathbf{e}
+\mathbf{b}^2\overline{\mathbf{a}})\mathbf{a}\mathbf{b})
-{\rm tr}(\overline{\mathbf{a}}\overline{\mathbf{b}}(t_{\overline{1}2^2}\overline{\mathbf{a}}\mathbf{b}^2-t_{1\overline{2}^2}\mathbf{e}
+\overline{\mathbf{b}}^2\mathbf{a}))  \\
=\ &r(t-\overline{t}^2)t_{121\overline{2}}+r(t\overline{t}-1)t_{1^22}+r\overline{t}(2\overline{s}-s^2)    \\
&+r(t^2-\overline{t})t_{\overline{1}\overline{2}\overline{1}2}+r(1-t\overline{t})t_{\overline{1}^2\overline{2}}+rt(\overline{s}^2-2s)   \\
&+t_{1\overline{2}^2}\big(\overline{t}t_{121\overline{2}}-tt_{1^22}+s^2-2\overline{s}\big)-t_{\overline{1}2^2}s+t_{2^3}  \\
&-t_{\overline{1}2^2}\big(tt_{\overline{1}\overline{2}\overline{1}2}-\overline{t}t_{\overline{1}^2\overline{2}}
+\overline{s}^2-2s\big)+t_{1\overline{2}^2}\overline{s}-t_{\overline{2}^3}   \\
=\ &\big(r(t-\overline{t}^2)+\overline{t}t_{1\overline{2}^2}\big)t_{121\overline{2}}
+\big(r(t^2-\overline{t})-tt_{\overline{1}2^2}\big)t_{\overline{1}\overline{2}\overline{1}2}  \\
&+(r(t\overline{t}-1)-tt_{1\overline{2}^2})t_{1^22}+(r(1-t\overline{t})+\overline{t}t_{\overline{1}2^2})t_{\overline{1}^2\overline{2}}  \\
&+(s^2-\overline{s})t_{1\overline{2}^2}+(s-\overline{s}^2)t_{\overline{1}2^2}
+r\overline{t}(2\overline{s}-s^2)+rt(\overline{s}^2-2s)+t^3-\overline{t}^3   \\
=\ &(rt-t^2\overline{t}+\overline{t}s)(t\overline{s}+(s+\overline{t})r+\overline{t}(1-t\overline{t}))  \\
&+(t\overline{t}^2-r\overline{t}-t\overline{s})(\overline{t}s+(t+\overline{s})r+t(1-t\overline{t}))+(t^3-r-ts)(ts-t\overline{t}+r)  \\
&+(r+\overline{t}\overline{s}-\overline{t}^3)(\overline{t}\overline{s}-t\overline{t}+r)
+(s^2-\overline{s})(\overline{t}r-t^2+s)  \\
&+(s-\overline{s}^2)(tr-\overline{t}^2+\overline{s})
+r\overline{t}(2\overline{s}-s^2)+rt(\overline{s}^2-2s)+t^3-\overline{t}^3   \\
=\ &s^3-\overline{s}^3+(r\overline{t}-2t^2)s^2-(rt-2\overline{t}^2)\overline{s}^2+(t^4+t^2\overline{t}+r^2t-r(t^2\overline{t}+3t))s  \\
&-(\overline{t}^4+\overline{t}^2t+r^2\overline{t}-r(\overline{t}^2t+3\overline{t}))\overline{s}+(t^3-\overline{t}^3)(r+1-t\overline{t}).
\end{align*}
In these computations, we have used the following identities:
\begin{alignat*}{2}
t_{2^3}&=t^3-3t\overline{t}+3, \qquad  &t_{\overline{2}^3}&=\overline{t}^3-3t\overline{t}+3,   \\
t_{1^22}&=ts-t\overline{t}+r,   \qquad   &t_{\overline{1}^2\overline{2}}&=\overline{t}\overline{s}-t\overline{t}+r,   \\
t_{\overline{1}2^2}&=tr-\overline{t}^2+\overline{s},  \qquad  &t_{1\overline{2}^2}&=\overline{t}r-t^2+s,    \\
t_{121\overline{2}}&=t\overline{s}+(s+\overline{t})r+\overline{t}(1-t\overline{t}),  \qquad
&t_{\overline{1}\overline{2}\overline{1}2}&=\overline{t}s+(t+\overline{s})r+t(1-t\overline{t}),
\end{alignat*}
which can be deduced by (\ref{eq:identity-1}), (\ref{eq:identity-2}).

\subsection{Proof of Theorem \ref{thm:main}}

Assume $t_{12\overline{12}}\ne t_{21\overline{21}}$, which implies that $(\mathbf{a},\mathbf{b})$ is irreducible.
Also assume that $\mathbf{a}\overline{\mathbf{b}}\overline{\mathbf{a}}\mathbf{b}\mathbf{a}
-\mathbf{b}\mathbf{a}\overline{\mathbf{b}}\overline{\mathbf{a}}\mathbf{b}$,
$\overline{\mathbf{b}}\mathbf{a}\mathbf{b}-\mathbf{a}\mathbf{b}\overline{\mathbf{a}}$, $\mathbf{a}-\mathbf{b}$ are not colinear, so $(\lambda,\mu,\nu)$ is determined up to a multiple. Furthermore, assume
$$\det(\lambda[\mathbf{a},\overline{\mathbf{b}}]+\nu\mathbf{e}+\nu[\overline{\mathbf{b}},\mathbf{a}])\ne0.$$
Then $\mathbf{y}$ is determined by the condition $\det(\mathbf{y})=1$ up to the multiplication by $e^{2k\pi i/3}$ with $0\le k\le 2$.

By results of \cite{La07}, the simultaneous conjugacy class of $(\mathbf{a},\mathbf{b})$ is determined by $t,\overline{t},s,\overline{s},r$ and 
$t_{12\overline{12}}$, subject to
$t_{12\overline{12}}^2-Pt_{12\overline{12}}+Q$ for some $P,Q\in\mathbb{C}[t,\overline{t},s,\overline{s},r]$. The assumption $t_{12\overline{12}}\ne t_{21\overline{21}}$ is equivalent to $P^2\ne 4Q$. 

Thus, a Zariski open subset of $\mathcal{X}_3^{\rm sym}(\Gamma)$ is a regular $\mathbb{Z}/6\mathbb{Z}$-cover of an Zariski open subset of the hypersurface in $\mathbb{C}^5$ defined by
\begin{align*}
s^3-\overline{s}^3+(r\overline{t}-2t^2)s^2-(rt-2\overline{t}^2)\overline{s}^2+(t^4+t^2\overline{t}+r^2t-r(t^2\overline{t}+3t))s   \\
-(\overline{t}^4+\overline{t}^2t+r^2\overline{t}-r(\overline{t}^2t+3\overline{t}))\overline{s}+(t^3-\overline{t}^3)(r+1-t\overline{t})=0.
\end{align*}

\bigskip

\noindent
{\bf Acknowledgement}.

I'd like to thank Professor Jiming Ma for helpful conversations. 

\bigskip

\noindent
Haimiao Chen (orcid: 0000-0001-8194-1264)\ \ \  \emph{chenhm@math.pku.edu.cn} \\
Department of Mathematics, Beijing Technology and Business University, \\
Liangxiang Higher Education Park, Fangshan District, Beijing, China.

\end{document}